\documentclass[12pt]{article}
\usepackage{amssymb,amsmath}
\usepackage{amsthm}
\usepackage{cases}
\usepackage{amsfonts}
\usepackage{graphicx}
\usepackage{cite,color,xcolor}
\usepackage[left=2.6cm,right=2.6cm,top=2.3cm,bottom=2.3cm]{geometry}
\usepackage[colorlinks,citecolor=blue,urlcolor=blue]{hyperref}
\usepackage[utf8]{inputenc}

\newtheorem{theorem}{Theorem}[section]

\newtheorem{lemma}{Lemma}[section]
\newtheorem{proposition}{Proposition}[section]

\newtheorem{definition}{Definition}[section]
\newtheorem{remark}{Remark}[section]

\usepackage{txfonts}
\newcommand{\bal}{\begin{align}}
\newcommand{\bbal}{\begin{align*}}
\newcommand{\beq}{\begin{equation}}
\newcommand{\eeq}{\end{equation}}
\newcommand{\bca}{\begin{cases}}
\newcommand{\eca}{\end{cases}}

\newcommand{\pa}{\partial}
\newcommand{\fr}{\frac}

\newcommand{\dd}{\mathrm{d}}

\newcommand{\R}{\mathbb{R}}

\newcommand{\Z}{\mathbb{Z}}

\newcommand{\les}{\lesssim}

\newcommand{\f}{\left}
\newcommand{\g}{\right}
\begin{document}
\bibliographystyle{plain}
\title{Zero-filter limit issue for the Camassa-Holm equation in Besov spaces}

\author{Yuxing Cheng$^{1}$, Jianzhong Lu$^{2,}$\footnote{E-mail: cyx191020@163.com; louisecan@163.com(Corresponding author); limin@jxufe.edu.cn; ny2008wx@163.com; lijinlu@gnnu.edu.cn
}, Min Li$^{3}$, Xing Wu$^{4}$, Jinlu Li$^{1}$\\
\small $^1$ School of Mathematics and Computer Sciences, Gannan Normal University, Ganzhou 341000, China\\
\small $^2$ School of Mathematics and Information, Xiangnan University, Chenzhou, 423000, China\\
\small $^3$ Department of Mathematics, Jiangxi University of Finance and Economics, Nanchang 330032, China
\\ \small $^4$ College of Information and Management Science, Henan Agricultural University, Zhengzhou 450002, China
}

\date{\today}
\maketitle\noindent{\hrulefill}

{\bf Abstract:} In this paper, we focus on zero-filter limit problem for the Camassa-Holm equation in the more general Besov spaces. We prove that the solution of the Camassa-Holm equation converges strongly in $L^\infty(0,T;B^s_{2,r}(\R))$ to the inviscid Burgers equation as the filter parameter $\alpha$ tends to zero with the given initial data $u_0\in B^s_{2,r}(\R)$. Moreover, we also show that the zero-filter limit for the Camassa-Holm equation does not converges uniformly with respect to the initial data in $B^s_{2,r}(\R)$.



\vskip0mm\noindent{\hrulefill}

\section{Introduction}

In this paper, we shall prove that the solution of
\begin{align}\label{alpha-c}
\begin{cases}
u_t+u \partial_xu=-\partial_x \left(1-\alpha^2 \partial_x^2\right)^{-1}\left(u^2+\frac{\alpha^2}{2} (\partial_xu)^2\right), \\
u(0,x)=u_0(x),
\end{cases}
\end{align}
which can be written as
\begin{align}\label{alpha-cr}
\begin{cases}
u_t+3u \partial_xu=-\alpha^2\partial^3_x \left(1-\alpha^2 \partial_x^2\right)^{-1}(u^2)-\frac{\alpha^2}{2}\partial_x \left(1-\alpha^2 \partial_x^2\right)^{-1} (\partial_xu)^2,\\
u(0,x)=u_0(x),
\end{cases}
\end{align}
converges to the solution of the following inviscid Burgers equation in the topology of Besov spaces $B^s_{2,r}(\R)$
\begin{align}\label{b}
\begin{cases}
u_t+3u \partial_xu=0,\\
u(0,x)=u_0(x).
\end{cases}
\end{align}

The equation \eqref{alpha-c} was firstly proposed in the context of hereditary symmetries studied in \cite{Fokas} as an integrable system. Then, it was rediscovered to describe the motion of shallow water waves by Camassa-Holm \cite{Camassa} and later became known as the Camassa-Holm(CH) equation. \eqref{alpha-c} is completely integrable \cite{Camassa,Constantin-P} with a bi-Hamiltonian structure \cite{Constantin-E,Fokas} and infinitely many conservation laws \cite{Camassa,Fokas}. Most importantly, CH equation has peakon solutions of the form $ce^{-|x-ct|}$ with $c>0,$ which aroused a lot of interest in physics, see \cite{Constantin-I,t}. Another remarkable feature of the CH equation is the wave breaking phenomena: the solution remains bounded while its slope becomes unbounded in finite time \cite{Constantin,Escher2,Escher3}. The existence of global weak solutions and dissipative solutions were investigated in \cite{bc1,bc2,xin}, more results can be found in the references therein.

Recently, a lot of literature has been devoted to studying the well-posedness problem of the Camassa–Holm equation in the Sobolev or Besov spaces. Li and Olver \cite{wp19,wp23} proved that the Cauchy problem (1.1) is locally well-posed with the initial data $u_0\in H^{s}(\R)$ and $s>3/2$. Danchin \cite{wp13} considered the local well-posedness in the Besov space, and proved the well-posedness in $B^s_{p,r}(\R)$ with $s>\max\{1+1/p,3/2\}$ and $1\leq p\leq\infty,~1\leq r<\infty$, while the continuous dependence part was proved by Li \cite{wp18}. For the endpoint $s=3/2$, local well-posedness in the space $B^{3/2}_{2,1}$ and ill-posedness in $B^{3/2}_{2,\infty}$ were also established by Danchin in \cite{wp14}. Finally, in the critical space $B^{3/2}_{2,2}=H^{3/2},$ Guo et al. \cite{glmy} proved a norm inflation and hence ill-posedness in $B^{1+1/p}_{p,r}(\R)$ with $1\leq p\leq\infty,~1< r\leq\infty$ by using peakon solutions, which solved the open problem left by Danchin \cite{wp14}.

Because of the interesting and remarkable features as mentioned above, the Camassa-Holm equation, as a class of integrable shallow water wave equations, has attracted much attention in recent twenty years. It is worth mentioning that the peakon solutions present the characteristic for the travelling water waves of greatest height and largest amplitude and arise as solutions to the free-boundary problem for incompressible Euler equations over a flat bed, see references \cite{Constantin-I,Escher4,Escher5} for the details.

Formally, in the limit case $\alpha=0,$ the CH equation becomes \eqref{b}, which can be seen as the non-viscous version ($\nu=0$) of the following well-known Burgers equation:
$$\partial_tu+u \partial_xu-\nu \partial_{xx}u=0.$$
This equation was proposed by Burgers in the 1940s as a toy model for turbulence. It mimics the Navier-Stokes equation of fluid motion through its fluid-like expressions for nonlinear advection and diffusion, but it is only one-dimensional and it lacks a pressure gradient driving the flow. For more details, please refer to \cite{miao2009,Molinet}. The Burgers equation has received extensive interests since 1940s, and it is perhaps the most basic example of a PDE evolution leading to shocks.

The zero-filter limit problem of the Camassa-Holm equation was first studied by Gui and Liu in \cite{GL}, where they proved that the solution of the dissipative Camassa-Holm equation does converge, at least locally, to the one of the dissipative Burgers equation as the filter parameter $\alpha$ tends to zero in the lower regularity Sobolev spaces. Very recently, Li et al. \cite{lyz}  proved that, as the filter parameter $\alpha\to 0,$ the solution of the Camassa-Holm equation with fractional dissipation converges strongly in $L^\infty(0,T;H^s(\R))$ to the inviscid Burgers equation for any given initial data $u_0\in H^s(\R)$ with $s>3/2$ and for some $T>0$. On the other hand, in \cite{lyz1}, the authors further showed that given initial data $u_0\in H^s(\R)$ with $s>3/2$ and for some $T>0$, the solution of the CH equation does not converges uniformly with respect to the initial data in $L^\infty(0,T;H^s(\R))$  to the inviscid Burgers equation as the filter parameter $\alpha$ tends to zero. In this paper, we mainly focus on the original Camassa-Holm equation \eqref{alpha-c}, and extend the previous two results to the more general Besov spaces.

Our first result proves that, given initial data $u_0\in B^s_{2,r}(\R)$ with $s>3/2, 1\leq r<\infty$ and for some $T>0$, the solution of the Camassa-Holm equation \eqref{alpha-c} converges strongly in $L^\infty(0,T;B^s_{2,r}(\R))$ to the inviscid Burgers equation \eqref{b} as the filter parameter $\alpha$ tends to zero. More precisely,
\begin{theorem}\label{th1} Let $s>\frac32, 1\leq r<\infty$ and $\alpha \in(0,1)$. Assume that the initial data $u_0\in B^s_{2,r}(\mathbb{R})$. Let $\mathbf{S}_{t}^{\mathbf{\alpha}}(u_0)$ and $\mathbf{S}_{t}^{0}(u_0)$ be the solutions of \eqref{alpha-c} and \eqref{b} with the initial data $u_0$, respectively. Then there exists a time $T=T(\|u_0\|_{B^s_{2,r}},s,r)>0$ such that  $\mathbf{S}_{t}^{\mathbf{\alpha}}(u_0),\mathbf{S}_{t}^{0}(u_0)\in \mathcal{C}([0,T];B^s_{2,r})$ and
$$
\lim_{\alpha\to 0}\left\|\mathbf{S}_{t}^{\mathbf{\alpha}}(u_0)
-\mathbf{S}_{t}^{0}(u_0)\right\|_{L^\infty_TB^{s}_{2,r}}=0.
$$
\end{theorem}

Our second result shows that, given initial data $u_0\in B^s_{2,r}(\R)$ with $s>3/2,1\leq r<\infty$ and for some $T>0$, the solution of the Camassa-Holm equation does not converges uniformly with respect to the initial data in $L^\infty(0,T;B^s_{2,r}(\R))$ to the inviscid Burgers equation as the filter parameter $\alpha$ tends to zero. To better state the result, we denote by $U_R$ the set of all functions $\phi$ in $B^s_{2,r}(\R)$ such that $\|\phi\|_{B^s_{2,r}(\R)}\leq R$, that is, 
$$U_R:=\big\{\phi\in B^s_{2,r}(\R): \|\phi\|_{B^s_{2,r}(\R)}\leq R\big\}.$$

\begin{theorem}\label{th2} Let $s>\frac32, 1\leq r<\infty$ and $\alpha\in(0,1)$. Let $\mathbf{S}_{t}^{\mathbf{\alpha}}(u_0)$ and $\mathbf{S}_{t}^{0}(u_0)$ be the solutions of \eqref{alpha-c} and \eqref{b} with the initial data $u_0\in B^s_{2,r}(\R)$, respectively. Then there exist a sequence of positive numbers $\{\alpha_n\}_{n=1}^{\infty}$ tending to zero and a sequence of initial data $\{u^n_0\}_{n=1}^{\infty}$ in $U_R$ such that for a short time $T_0\leq T$
$$\liminf_{\alpha_n\to0}\left\|\mathbf{S}_{t}^{\alpha_n}(u^n_0)-\mathbf{S}_{t}^{0}(u^n_0)\right\|_{
L^\infty_{T_0}B^s_{2,r}}\geq \eta_0,$$
with some positive constant $\eta_0$.
\end{theorem}

\begin{remark}
Note that if  $r=2,$ then $H^s(\R)= B^s_{2,2}$. Hence, our results improve the corresponding results in \cite{lyz} and \cite{lyz1}.
\end{remark}

\section{Preliminaries}\label{sec2}
{\bf Notation}\; Throughout this paper, we denote by $C$ any positive constant independent of the parameter $\alpha$, which may change from line to line.
Given a Banach space $X$, we denote its norm by $\|\cdot\|_{X}$. For $I\subset\R$, we denote by $\mathcal{C}(I;X)$ the set of continuous functions on $I$ with values in $X$. Sometimes, we also denote $L^p(0,T;X)$ by $L_T^pX$.
For all $f\in \mathcal{S}'$, the Fourier transform  $\mathcal{F}(f)$ and inverse Fourier transform  $\mathcal{F}^{-1}(f)$ are defined respectively by
$$
\mathcal{F}(f)(\xi)=\int_{\R}e^{-ix\xi}f(x)\dd x, \quad \mathcal{F}^{-1}(f)(\xi)=\int_{\R}e^{-ix\xi}f(x)\dd x, \qquad\text{for any}\; \xi\in\R.
$$
Next, we recall some facts about the Littlewood-Paley decomposition, the nonhomogeneous Besov spaces and several useful related properties.
\begin{proposition}[Littlewood-Paley decomposition, See \cite{B.C.D}] Let $\mathcal{B}:=\{\xi\in\mathbb{R}:|\xi|\leq \frac 4 3\}$ and $\mathcal{C}:=\{\xi\in\mathbb{R}:\frac 3 4\leq|\xi|\leq \frac 8 3\}.$
There exist two radial functions $\chi\in C_c^{\infty}(\mathcal{B})$ and $\varphi\in C_c^{\infty}(\mathcal{C})$ both taking values in $[0,1]$ such that
\begin{align*}
&\chi(\xi)+\sum_{j\geq0}\varphi(2^{-j}\xi)=1 \quad \forall \;  \xi\in \R,\\
&\frac{1}{2} \leq \chi^{2}(\xi)+\sum_{j \geq 0} \varphi^{2}(2^{-j} \xi) \leq 1\quad \forall \;  \xi\in \R.
\end{align*}
\end{proposition}

\begin{definition}[See \cite{B.C.D}]
For every $u\in \mathcal{S'}(\mathbb{R})$, the Littlewood-Paley dyadic blocks ${\Delta}_j$ are defined as follows
\begin{numcases}{\Delta_ju=}
0, & if $j\leq-2$;\nonumber\\
\chi(D)u=\mathcal{F}^{-1}(\chi \mathcal{F}u), & if $j=-1$;\nonumber\\
\varphi(2^{-j}D)u=\mathcal{F}^{-1}\big(\varphi(2^{-j}\cdot)\mathcal{F}u\big), & if $j\geq0$.\nonumber
\end{numcases}
In the inhomogeneous case, the following Littlewood-Paley decomposition makes sense
$$
u=\sum_{j\geq-1}{\Delta}_ju,\quad \forall\;u\in \mathcal{S'}(\mathbb{R}).
$$
\end{definition}
\begin{definition}[See \cite{B.C.D}]
Let $s\in\mathbb{R}$ and $(p,q)\in[1, \infty]^2$. The nonhomogeneous Besov space $B^{s}_{p,q}(\R)$ is defined by
\begin{align*}
B^{s}_{p,q}(\R):=\Big\{f\in \mathcal{S}'(\R):\;\|f\|_{B^{s}_{p,q}(\mathbb{R})}<\infty\Big\},
\end{align*}
where
\begin{numcases}{\|f\|_{B^{s}_{p,q}(\mathbb{R})}=}
\left(\sum_{j\geq-1}2^{sjq}\|\Delta_jf\|^r_{L^p(\mathbb{R})}\right)^{\fr1q}, &if $1\leq q<\infty$,\nonumber\\
\sup_{j\geq-1}2^{sj}\|\Delta_jf\|_{L^p(\mathbb{R})}, &if $q=\infty$.\nonumber
\end{numcases}
\end{definition}
\begin{remark}\label{re3}
It should be emphasized that the following embedding will be often used implicity:
$$B^s_{p,q}(\R)\hookrightarrow B^t_{p,r}(\R)\quad\text{for}\;s>t\quad\text{or}\quad s=t,1\leq q\leq r\leq\infty.$$
\end{remark}
\begin{lemma}[\cite{B.C.D}]\label{le-pro}
For $s>0$ and $1\leq q\leq \infty$, $B^s_{2,q}(\R)\cap L^\infty(\R)$ is an algebra.
Moreover, we have for any $u,v \in B^s_{2,q}(\R)\cap L^\infty(\R)$
\begin{align*}
&\|uv\|_{B^s_{2,q}(\R)}\leq C\big(\|u\|_{B^s_{2,q}(\R)}\|v\|_{L^\infty(\R)}+\|v\|_{B^s_{2,q}(\R)}\|u\|_{L^\infty(\R)}\big).
\end{align*}
In particular, for $s>\frac12$, due to the fact $B^s_{2,q}(\R)\hookrightarrow L^\infty(\R)$, then we have
\begin{align*}
&\|uv\|_{B^s_{2,q}(\R)}\leq C\|u\|_{B^s_{2,q}(\R)}\|v\|_{B^s_{2,q}(\R)}.
\end{align*}
\end{lemma}
\begin{lemma}[\cite{B.C.D}]\label{le-pro1}
For $s>\frac32$ and $1\leq q\leq \infty$, there holds
\begin{align*}
&\|uv\|_{B^{s-2}_{2,q}(\R)}\leq C\|u\|_{B^{s-2}_{2,q}(\R)}\|v\|_{B^{s-1}_{2,q}(\R)}.
\end{align*}
\end{lemma}
\begin{lemma}[\cite{B.C.D}]\label{le3}
Let $s>0,1\leq q\leq \infty$ and $f,g\in {\rm{Lip}(\R)}\cap B^s_{2,q}(\R)$. Then we have
\bbal
\big|\big|(2^{js}||[\Delta_j,f]\pa_xg||_{L^2(\mathbb{R})})_{j\geq-1}\big|\big|_{\ell^q}\leq
\begin{cases}
C\big(\|\partial_x f\|_{L^\infty(\R)}\|g\|_{B^{s}_{2,q}(\R)}+\|f\|_{B^s_{2,q}(\R)}\|\pa_xg\|_{L^\infty(\R)}\big),\\
C\big(\|\partial_x f\|_{L^\infty(\R)}\|g\|_{B^{s}_{2,q}(\R)}+\|\pa_x f\|_{B^s_{2,q}(\R)}\|g\|_{L^\infty(\R)}\big).
\end{cases}
\end{align*}
\end{lemma}

\begin{lemma}[\cite{Li2020}]\label{le2} Let $s\in\R$ and $1\leq q\leq \infty$. Define the high frequency function $f_n$ and the low frequency function $g_n$ by
\bbal
&f_n=2^{-ns}\phi(x)\sin \f(\frac{17}{12}2^nx\g), \qquad g_n=2^{-n}\phi(x),\quad n\gg1.
\end{align*}
where $\widehat{\phi}\in \mathcal{C}^\infty_0(\mathbb{R})$ be an even, real-valued and non-negative function on $\R$ and satisfy
\begin{numcases}{\widehat{\phi}(\xi)=}
1,&if $|\xi|\leq \frac{1}{4}$,\nonumber\\
0,&if $|\xi|\geq \frac{1}{2}$.\nonumber
\end{numcases}
Then for any $\sigma\in\R$, there exists two positive constants $c=c(\phi)$ and $C=C(\phi)$ such that
\bal
&\|f_n\|_{L^\infty(\R)}\leq C2^{-ns},\quad\|\pa_xf_n\|_{L^\infty(\R)}\leq C2^{-n(s-1)},\label{m2}\\
&\|f_n\|_{B_{2,q}^\sigma(\R)}\approx2^{n(\sigma-s)},\quad \|g_n\|_{B_{2,q}^\sigma(\R)}\approx2^{-n},\label{m3}\\
&\liminf_{n\rightarrow \infty}\|g_n\pa_xf_n\|_{B^s_{2,q}(\R)}\geq c. \label{m6}
\end{align}
\end{lemma}

\begin{lemma}\label{lemm1}
For any $\alpha \in(0,1)$, there holds
\bbal
\Big|\int_{\R}(1-\alpha^2\pa^2_x)^{-1}\Delta_j(vu_x)\Delta_ju\dd x\Big|\leq C\|\pa_xv\|_{L^\infty(\R)}\|\Delta_ju\|^2_{L^2(\R)}+C\|[\Delta_j,v]\pa_xu\|_{L^2(\R)}\|\Delta_ju\|_{L^2(\R)}.
\end{align*}
\end{lemma}
\begin{proof}
To bound this term, we can rewrite it as
\begin{align}\label{es-n1}
\int_{\R}(1-\alpha^2\pa^2_x)^{-1}\Delta_j(vu_x)\Delta_ju\dd x&=\int_{\R}(1-\alpha^2\pa^2_x)^{-1}( v\pa_x\Delta_ju)\cdot \Delta_ju\dd x \nonumber\\& \quad +\int_{\R}(1-\alpha^2\pa^2_x)^{-1}[\Delta_j,v]\pa_xu\cdot \Delta_ju\dd x:=\mathrm{I}+\mathrm{II}.
\end{align}
Letting $u_j=(1-\alpha^2\pa^2_x)^{-1}\Delta_ju$, then we have
\bbal
\mathrm{I}&=\int_{\R}( v\pa_x\Delta_ju)\cdot (1-\alpha^2\pa^2_x)^{-1}\Delta_ju\dd x
\\&=\int_{\R} v(1-\alpha^2\pa^2_x)\pa_xu_j\cdot u_j\dd x
\\&=-\frac12\int_{\R} \pa_xv\cdot u_j^2\dd x+\int_{\R} \pa_xv(\alpha^2\pa^2_xu_j)\cdot u_j\dd x-\fr12\int_{\R} \pa_xv\cdot (\alpha\pa_xu_j)^2\dd x,
\end{align*}
which along with Holder's inequality yields
\bal\label{es-i}
|\mathrm{I}|&\leq C\|\pa_xv\|_{L^\infty}\f(\|u_j\|^2_{L^2}+\|\alpha^2\pa^2_xu_j\|^2_{L^2}+\|\alpha\pa_xu_j\|^2_{L^2}\g)\\
& \nonumber \leq C\|\pa_xv\|_{L^\infty}\|\Delta_ju\|^2_{L^2}.
\end{align}
According to Holer's inequality, we also have
\bal\label{es-ii}
|\mathrm{II}|&=|\int_{\R} [\Delta_j,v]\pa_xu\cdot u_j\dd x|
\leq \|[\Delta_j,v]\pa_xu\|_{L^2}\|u_j\|_{L^2}
\leq C\|[\Delta_j,v]\pa_xu\|_{L^2}\|\Delta_ju\|_{L^2}.
\end{align}
Inserting \eqref{es-i} and \eqref{es-ii} into \eqref{es-n1}, we completes the proof of this lemma.
\end{proof}

\begin{lemma}\label{lemm2}
For any $\alpha \in(0,1)$, there holds
\bbal
\Big|\frac{\alpha^2}{2}\int_{\R}\pa_x(1-\alpha^2\pa^2_x)^{-1}\Delta_j(\partial_xu\pa_xv)\cdot \Delta_jw\dd x\Big|\leq C2^{-j}\|\Delta_j(\partial_xu\pa_xv)\|_{L^2}\|\Delta_jw\|_{L^2}.
\end{align*}
\end{lemma}
\begin{proof}
It is easy to show that
\bbal
\Big|\frac{\alpha^2}{2}\int_{\R}\pa_x(1-\alpha^2\pa^2_x)^{-1}\Delta_j(\partial_xu\pa_xv)\cdot \Delta_jw\dd x\Big|\leq \|\Delta_jw\|_{L^2}\alpha^2||\pa_x(1-\alpha^2\pa^2_x)^{-1}\Delta_j(\partial_xu\pa_xv)||_{L^2}.
\end{align*}
By the definition of $\Delta_j$ and Plancherel's identity, we have
\bbal
||\pa_x(1-\alpha^2\pa^2_x)^{-1}\Delta_j(\partial_xu\pa_xv)||_{L^2}&=||\frac{i\xi}{1+\alpha^2\xi^2}\mathcal{F}[\Delta_j(\partial_xu\pa_xv)]||_{L^2}
\\&\leq C\frac{1}{\alpha^2}2^{-j}||\mathcal{F}[\Delta_j(\partial_xu\pa_xv)]||_{L^2}
\\&\leq \frac{1}{\alpha^2}2^{-j}||\Delta_j(\partial_xu\pa_xv)||_{L^2}.
\end{align*}
Combining this two estimates, we complete the proof of this lemma.
\end{proof}

\begin{lemma}\label{lemm3}
Let $\sigma\in \R$ and $r\in [1,\infty)$. For any $\alpha \in(0,1)$, there holds
\bbal
&||\alpha^2\pa^2_x(1-\alpha^2\pa^2_x)^{-1}u||_{B^\sigma_{2,r}}+||\alpha\pa_x(1-\alpha^2\pa^2_x)^{-1}u||_{B^\sigma_{2,r}}\leq C||u||_{B^\sigma_{2,r}},
\\& ||\alpha^2\pa_x(1-\alpha^2\pa^2_x)^{-1}u||_{B^\sigma_{2,r}}\leq C||u||_{B^{\sigma-1}_{2,r}}.
\end{align*}
\end{lemma}
\begin{proof}
The results can easily deduce from the definition of $\Delta_j$ and Plancherel's identity. Here, we omit it.
\end{proof}

\section{Proof of Theorem \ref{th1}}
 \quad In this section, we will prove that the zero-filter limit issue for the Camassa-Holm equation converges strongly in $L^\infty(0,T;B^s_{2,r}(\R))$ to the inviscid Burgers equation. We divide the proof of Theorem \ref{th1} into three steps.

 {\bf Step 1: Uniform bound w.r.s $\alpha\in(0,1)$ of $\mathbf{S}_{t}^{\mathbf{\alpha}}(u_0)$ in $B^s_{2,r}$}. For fixed $\alpha>0$, by the classical local well-posedness result, we known that there exists a $T_\alpha=T(\|u_0\|_{B^s_{2,r}},s,r,\alpha)>0$ such that the Camassa-Holm has a unique solution $\mathbf{S}_{t}^{\mathbf{\alpha}}(u_0)\in\mathcal{C}([0,T_\alpha];B^s_{2,r})$.

We shall prove that $\exists\; T=T(\|u_0\|_{B^s_{2,r}},s,r)>0$ such that $T\leq T_{\alpha}$ and there exists $C>0$ independent of $\alpha$ such that
\begin{align}\label{m1}
\|\mathbf{S}_{t}^{\mathbf{\alpha}}(u_0)\|_{L_T^{\infty} B^s_{2,r}} \leq C, \quad \forall \alpha \in[0,1).
\end{align}
To simplify notation, we set $u=\mathbf{S}_{t}^{\mathbf{\alpha}}(u_0)$.
Applying the operator $\Delta_j$ to $\eqref{alpha-c}$, multiplying $\Delta_ju$  and integrating the resulting over $\R$, we obtain
\begin{align}
\frac12\frac{\dd }{\dd t}\|\Delta_ju\|^2_{L^2}&=\frac12\int_{\R}\pa_xu|\Delta_ju|^2\dd x-\int_{\R}[\Delta_j,u]\pa_xu\cdot \Delta_ju\dd x\label{y1}\\
&\quad-2\int_{\R}(1-\alpha^2\pa^2_x)^{-1}\Delta_j(uu_x)\cdot \Delta_ju\dd x\label{y2}\\
&\quad
-\frac{\alpha^2}{2}\int_{\R}\pa_x(1-\alpha^2\pa^2_x)^{-1}\Delta_j(\partial_xu)^2\cdot \Delta_ju\dd x.\label{y3}
\end{align}
To bound \eqref{y1}, it is easy to obtain
\bbal
|\eqref{y1}|
\leq C\|\pa_xu\|_{L^\infty}\|\Delta_ju\|^2_{L^2}+C\|[\Delta_j,u]\pa_x u\|_{L^2}\|\Delta_ju\|_{L^2},
\end{align*}
To bound \eqref{y2}, by Lemma \ref{lemm1}, we have
\bbal
|\eqref{y2}|\leq C\|\pa_xu\|_{L^\infty}\|\Delta_ju\|^2_{L^2}+C\|[\Delta_j,u]\pa_xu\|_{L^2}\|\Delta_ju\|_{L^2}
\end{align*}
To bound \eqref{y3}, we can deduce from Lemma \ref{lemm2} that
\bbal
|\eqref{y3}|\leq C2^{-j}\|\Delta_j(\partial_xu)^2\|_{L^2}\|\Delta_ju\|_{L^2}.
\end{align*}
Combining the above yields that
\begin{align*}
\frac{\dd }{\dd t}\|\Delta_ju\|_{L^2}\leq C\Big(\|\pa_xu\|_{L^\infty}\|\Delta_ju\|_{L^2}+\|[\Delta_j,u]\pa_xu\|_{L^2}
+2^{-j}\|\Delta_j(\partial_xu)^2\|_{L^2}\Big)
\end{align*}
which leads to
\bbal
\|\Delta_ju\|_{L^2}\leq \|\Delta_ju_0\|_{L^2}+\int^t_0\Big(\|\pa_xu\|_{L^\infty}\|\Delta_ju\|_{L^2}+\|[\Delta_j,u]\pa_xu\|_{L^2}
+2^{-j}\|\Delta_j(\partial_xu)^2\|_{L^2}\Big)\dd \tau.
\end{align*}
Multiplying the above inequality  by $2^{js}$ and taking the $\ell^r$ norm over $\Z$, we obtain from Lemma \ref{le-pro} that
\bbal
\|u(t)\|_{B^s_{2,r}}&\leq \|(u_0)\|_{B^s_{2,r}} + C\int^t_0\Big(\|\pa_xu\|_{L^\infty}\|u\|_{B^s_{2,r}}+\|(\pa_xu)^2\|_{B^{s-1}_{2,r}}\Big)\dd \tau
\\&\leq \|(u_0)\|_{B^s_{2,r}} + C\int^t_0\|u\|^2_{B^s_{2,r}}\dd \tau.
\end{align*}
Thus, by continuity arguments there exists a time $T=T(\|u_0\|_{B^s_{2,r}})>0$ such that \eqref{m1} holds uniformly w.r.s $\alpha\in(0,1)$. Moreover, if $u_0 \in H^\gamma$ for some $\gamma>s$, then there exists $C_2>0$ independent of $\alpha$ such that
\begin{align}\label{es-u-hi}
\|\mathbf{S}_{t}^{\mathbf{\alpha}}(u_0)\|_{L_T^{\infty} H^\gamma} \leq C_2\left\|u_0\right\|_{H^\gamma} .
\end{align}

{\bf Step 2: Estimations of $\|\mathbf{S}_{t}^{\mathbf{\alpha}}(u_0)-\mathbf{S}_{t}^{\alpha}(S_nu_0)\|_{B^s_{2,r}}$ .}  Denoting $$\mathbf{v}(t)=\mathbf{S}_{t}^{\mathbf{\alpha}}(u_0)-\mathbf{S}_{t}^{\alpha}(S_nu_0)\quad\text{and}\quad \mathbf{v}|_{t=0}=(\mathrm{Id}-S_n)u_0,$$ we infer that $\mathbf{v}$ satisfies
$$
\partial_{t} \mathbf{v}+\mathbf{S}_{t}^{\mathbf{\alpha}}(u_0)\partial_{x} \mathbf{v}=-\mathbf{v}\partial_{x} \mathbf{S}_{t}^{\alpha}(S_nu_0) +\mathcal{B}(\mathbf{v}, \mathbf{S}_{t}^{\mathbf{\alpha}}(u_0)+\mathbf{S}_{t}^{\alpha}(S_nu_0)),
$$
where
$$
\mathcal{B}:(f, g) \mapsto \partial_x\left(1-\alpha^2 \partial_x^2\right)^{-1}\left(f g+\frac{\alpha^2}{2} \partial_{x} f \partial_{x} g\right).
$$
Notice that
\bal
&\quad \ \int_{\R}\Delta_j\f(\mathcal{B}(\mathbf{v}, \mathbf{S}_{t}^{\mathbf{\alpha}}(u_0)+\mathbf{S}_{t}^{\alpha}(S_nu_0))\g)\cdot \Delta_j\mathbf{v}\dd x \nonumber
\\&=\frac12\int_{\R}\Delta_j\f(\alpha^2\pa_x(1-\alpha^2\pa^2_x)^{-1}(\mathbf{v}_x \pa_x[\mathbf{S}_{t}^{\mathbf{\alpha}}(u_0)+\mathbf{S}_{t}^{\alpha}(S_nu_0)])\g)\cdot \Delta_j\mathbf{v}\dd x \label{lyz+1}
\\&\quad +2\int_{\R}(1-\alpha^2\pa^2_x)^{-1} \Delta_j(\pa_x\mathbf{S}_{t}^{\alpha}(S_nu_0)\mathbf{v})\cdot \Delta_j\mathbf{v}\dd x\label{lyz+2}\\
&\quad +2\int_{\R}(1-\alpha^2\pa^2_x)^{-1}\Delta_j
\f[\f(\mathbf{v}+\mathbf{S}_{t}^{\alpha}(S_nu_0)\g)\mathbf{v}_x)\g]\cdot \Delta_j\mathbf{v}\dd x. \label{lyz+3}
\end{align}
To bound \eqref{lyz+1} and \eqref{lyz+3}, we easily from Lemmas \ref{lemm1}-\ref{lemm2} that
\bbal
|\eqref{lyz+1}|&\leq C2^{-j}||\Delta_j(\mathbf{v}_x \pa_x[\mathbf{S}_{t}^{\mathbf{\alpha}}(u_0)+\mathbf{S}_{t}^{\alpha}(S_nu_0)])||_{L^2}\|\Delta_j\mathbf{v}\|_{L^2},\\
|\eqref{lyz+3}|&\leq C\|\pa_x[\mathbf{v}+\mathbf{S}_{t}^{\alpha}(S_nu_0)]\|_{L^\infty}\|\Delta_j\mathbf{v}\|^2_{L^2}
+C\|[\Delta_j,\mathbf{v}+\mathbf{S}_{t}^{\alpha}(S_nu_0)]\pa_x\mathbf{v}\|_{L^2}\|\Delta_j\mathbf{v}\|_{L^2}
\end{align*}
By Holer's inequality, we have
\bbal
|\eqref{lyz+2}|\leq C||\Delta_j(\pa_x\mathbf{S}_{t}^{\alpha}(S_nu_0)\mathbf{v})||_{L^2}|| \Delta_j\mathbf{v}||_{L^2}.
\end{align*}
Following the same procedure as that in Step 1, we deduce
\bal\label{chen1}
&\|\Delta_j\mathbf{v}\|_{L^2}\leq \|(Id-S_n)u_0\|_{L^2}+\int^t_0\Big(\|\pa_x\mathbf{S}_{t}^{\mathbf{\alpha}}(u_0)\|_{L^\infty}\|\Delta_j\mathbf{v}\|_{L^2}+\|[\Delta_j,\mathbf{S}_{t}^{\mathbf{\alpha}}(u_0)]\pa_x\mathbf{v}\|_{L^2}
\\&\nonumber\qquad +||\Delta_j(\pa_x\mathbf{S}_{t}^{\alpha}(S_nu_0)\mathbf{v})||_{L^2}+\|[\Delta_j,\mathbf{v}+\mathbf{S}_{t}^{\mathbf{\alpha}}(u_0)]\pa_x\mathbf{v}\|_{L^2}
+2^{-j}\|\Delta_j(\mathbf{v}_x \pa_x[\mathbf{S}_{t}^{\mathbf{\alpha}}(u_0)+\mathbf{S}_{t}^{\alpha}(S_nu_0)])\|_{L^2}\Big)\dd \tau.
\end{align}
Multiplying \eqref{chen1}  by $2^{js}$ and taking the $\ell^r$ norm over $\Z$, we obtain from Lemma \ref{le-pro} that
\bal\label{chen-vs}
\|\mathbf{v}(t)\|_{B^s_{2,r}}&\leq \|(Id-S_n)u_0\|_{B^s_{2,r}} + C\int^t_0\Big(\|\mathbf{v},\mathbf{S}_{t}^{\mathbf{\alpha}}(u_0)\|_{B^s_{2,r}}\|\mathbf{v}\|_{B^s_{2,r}}+\|\mathbf{v}_x \pa_x[\mathbf{S}_{t}^{\mathbf{\alpha}}(u_0)+\mathbf{S}_{t}^{\alpha}(S_nu_0)]\|_{B^{s-1}_{2,r}}\Big)\dd \tau \nonumber
\\& \qquad +\int^t_0||\pa_x\mathbf{S}_{t}^{\alpha}(S_nu_0)\mathbf{v}||_{B^s_{2,r}}\dd \tau \nonumber
\\&\leq \|(Id-S_n)u_0\|_{B^s_{2,r}} + C\int^t_0\|\mathbf{v}\|_{B^s_{2,r}}\dd \tau+C2^n\int^t_0\|\mathbf{v}\|_{B^{s-1}_{2,r}}\dd \tau.
\end{align}
The above inequality reduces to
\begin{align}\label{chen-vsr}
\|\mathbf{v}\|_{{B}^{s}_{2,r}}
&\leq \|(Id-S_n)u_0\|_{B^s_{2,r}}+ C2^n\int^t_0\|\mathbf{v}\|_{B_{2,r}^{s-1}}\dd \tau.
\end{align}
To close \eqref{chen-vsr}, we have to estimate $\|\mathbf{v}\|_{B_{2,r}^{s-1}}$. Similar argument as \eqref{chen-vs}, we have from Lemma \ref{le-pro1} that
\bbal
\|\mathbf{v}(t)\|_{B^{s-1}_{2,r}}&\leq \|(Id-S_n)u_0\|_{B^{s-1}_{2,r}} + C\int^t_0\Big(\|\mathbf{v},\mathbf{S}_{t}^{\mathbf{\alpha}}(u_0)\|_{B^{s}_{2,r}}\|\mathbf{v}\|_{B^{s-1}_{2,r}}+\|\mathbf{v}_x \pa_x[\mathbf{S}_{t}^{\mathbf{\alpha}}(u_0)+\mathbf{S}_{t}^{\alpha}(S_nu_0)]\|_{B^{s-2}_{2,r}}\Big)\dd \tau
\\& \qquad +\int^t_0||\pa_x\mathbf{S}_{t}^{\alpha}(S_nu_0)\mathbf{v}||_{B^{s-1}_{2,r}}\dd \tau
\\&\leq \|(Id-S_n)u_0\|_{B^s_{2,r}} + C\int^t_0\|\mathbf{v}\|_{B^{s-1}_{2,r}}\dd \tau.
\end{align*}
Applying Gronwall's inequality yields that for $t\in [0,T]$
\bal\label{l4}
\|\mathbf{v}(t)\|_{B_{2,r}^{s-1}}\leq C\|(\mathrm{Id}-S_n)u_0\|_{B_{2,r}^{s-1}}\leq C2^{-n}\|(\mathrm{Id}-S_n)u_0\|_{B_{2,r}^{s}}.
\end{align}
Inserting \eqref{l4} into \eqref{chen-vsr} and applying Gronwall's inequality, we obtain that for $t\in [0,T]$
\bal\label{solu-diff}
\|\mathbf{S}_{t}^{\mathbf{\alpha}}(u_0)-\mathbf{S}_{t}^{\alpha}(S_nu_0)\|_{B_{2,r}^{s}}\leq C\|(\mathrm{Id}-S_n)u_0\|_{B_{2,r}^{s}}, \quad \forall \alpha\in[0,1).
\end{align}

{\bf Step 3: Estimation of $\|\mathbf{S}_{t}^{\alpha}(S_nu_0)-\mathbf{S}_{t}^{0}(S_nu_0)\|_{B^s_{2,r}}$}. We can find that $\mathbf{S}_{t}^{\alpha}(S_nu_0)$ satisfies the following equation
\bbal
\pa_t\mathbf{S}_{t}^{\alpha}(S_nu_0)+3\mathbf{S}_{t}^{\alpha}(S_nu_0)\pa_x\mathbf{S}_{t}^{\alpha}(S_nu_0)
&=-\frac12\alpha^2\pa_x(1-\alpha^2\pa^2_x)^{-1}[\pa_x\mathbf{S}_{t}^{\alpha}(S_nu_0)]^2
\\&\quad -\alpha^2\pa^3_x(1-\alpha^2\pa^2_x)^{-1}[\mathbf{S}_{t}^{\alpha}(S_nu_0)]^2.
\end{align*}
Denoting $$\mathbf{w}(t)= \mathbf{S}_{t}^{\alpha}(S_nu_0)-\mathbf{S}_{t}^{0}(S_nu_0)\quad\text{and}\quad \mathbf{w}|_{t=0}=0,$$ we infer that $\mathbf{w}$ satisfies
$$
\left\{\begin{array}{l}
\partial_t\mathbf{w}+3\mathbf{S}_{t}^{0}(S_nu_0) \partial_x \mathbf{w}=-3\mathbf{w} \partial_x \mathbf{S}_{t}^{\alpha}(S_nu_0)-\mathbf{I}, \\
\mathbf{w}|_{t=0}=0,
\end{array}\right.
$$
where
\bbal
\mathbf{I}:=\alpha^2\partial_x\left(1-\alpha^2 \partial_x^2\right)^{-1}\f[\partial^2_x \left([\mathbf{S}_{t}^{\alpha}(S_nu_0)]^2\right)+\frac{1}{2}\left[\pa_x\mathbf{S}_{t}^{\alpha}(S_nu_0)\right]^2\g].
\end{align*}
For the $\mathbf{I}$, it is easy to show that for $j\geq -1$,
\bal\label{es-I}
\Big|\int_{\R}\Delta_j\mathbf{I}\cdot \Delta_j\mathbf{w}\dd x\Big|\leq \alpha\|\Delta_j\mathbf{w}\|_{L^2}||\Delta_j\big(\partial^2_x \left([\mathbf{S}_{t}^{\alpha}(S_nu_0)]^2\right)+\frac{1}{2}\left[\pa_x\mathbf{S}_{t}^{\alpha}(S_nu_0)\right]^2\big)||_{L^2}.
\end{align}
Taking the similar argument with Step 1 and using \eqref{es-I}, we have
\bbal
\|\Delta_j\mathbf{w}\|_{L^2}&\leq C\int^t_0\Big(\|\pa_x\mathbf{S}_{t}^{0}(S_nu_0)\|_{L^\infty}\|\Delta_j\mathbf{w}\|_{L^2}
+\|[\Delta_j,\mathbf{S}_{t}^{0}(S_nu_0)]\pa_x\mathbf{w}\|_{L^2}
\\&\qquad +\|\Delta_j[\mathbf{w} \partial_x \mathbf{S}_{t}^{\alpha}(S_nu_0)]\|_{L^2}
+\alpha\|\Delta_j\big(\partial^2_x \left([\mathbf{S}_{t}^{\alpha}(S_nu_0)]^2\right)+\frac{1}{2}\left[\pa_x\mathbf{S}_{t}^{\alpha}(S_nu_0)\right]^2\big)\|_{L^2}\Big)\dd \tau.
\end{align*}
Multiplying the above inequlity  by $2^{j(s-1)}$ and taking the $\ell^r$ norm over $\Z$, we obtain from Lemma \ref{le-pro} that
\bal\label{es-w}
\|\mathbf{w}(t)\|_{B^{s-1}_{2,r}}&\leq  C\int^t_0\Big(\|\mathbf{S}_{t}^{\mathbf{0}}(S_nu_0)\|_{B^{s}_{2,r}}\|\mathbf{w}\|_{B^{s-1}_{2,r}}
+\alpha\|\partial^2_x \left([\mathbf{S}_{t}^{\alpha}(S_nu_0)]^2\right)+\frac{1}{2}\left[\pa_x\mathbf{S}_{t}^{\alpha}(S_nu_0)\right]^2\|_{B^{s-1}_{2,r}}\Big)\dd \tau\nonumber
\\& \qquad +\int^t_0||\pa_x\mathbf{S}_{t}^{\alpha}(S_nu_0)\mathbf{w}||_{B^{s-1}_{2,r}}\dd \tau\nonumber
\\&\leq C\int^t_0\|\mathbf{w}\|_{B^{s-1}_{2,r}}\dd \tau+C\alpha\int^t_0\|\mathbf{S}_{t}^{\alpha}(S_nu_0),\pa_x\mathbf{S}_{t}^{\alpha}(S_nu_0)\|_{L^\infty}\|\mathbf{S}_{t}^{\alpha}(S_nu_0)\|_{B^{s+1}_{2,r}}\dd \tau \nonumber
\\&\leq C\int^t_0\|\mathbf{w}\|_{B^{s-1}_{2,r}}\dd \tau+C\alpha\int^t_0\|\mathbf{S}_{t}^{\alpha}(S_nu_0)\|_{B^{s}_{2,r}}\|\mathbf{S}_{t}^{\alpha}(S_nu_0)\|_{B^{s+1}_{2,r}}\dd \tau.
\end{align}
Combining \eqref{es-w} and \eqref{es-u-hi},  we deduce that
\begin{align*}
\|\mathbf{w}(t)\|_{B_{2,r}^{s-1}}
&\leq C \int_{0}^{t}\|\mathbf{w}\|_{B_{2,r}^{s-1}} \dd\tau+C\alpha2^{n},
\end{align*}
which along with Gronwall's inequality implies
\bbal
\|\mathbf{w}(t)\|_{B_{2,r}^{s-1}}\leq C\alpha2^{n}.
\end{align*}
Then, we get for $t\in[0,T]$
\begin{align}\label{es-w1}
\|\mathbf{w}(t)\|_{B_{2,r}^{s}} &\leq \|\mathbf{w}(t)\|^{\frac12}_{B_{2,r}^{s-1}}\|\mathbf{w}(t)\|^{\frac12}_{B_{2,r}^{s+1}} \leq C \alpha^{\fr12} 2^{n}.
\end{align}
Due to \eqref{es-w1} and \eqref{solu-diff}, we have for $t\in[0,T]$
\bbal
&\|\mathbf{S}_{t}^{\mathbf{\alpha}}(u_0)-\mathbf{S}_{t}^{0}(u_0)\|_{B_{2,r}^{s}}\\
&\leq \|\mathbf{S}_{t}^{\mathbf{\alpha}}(u_0)-\mathbf{S}_{t}^{\alpha}(S_nu_0)\|_{B_{2,r}^{s}}+\|\mathbf{S}_{t}^{\alpha}(S_nu_0)-\mathbf{S}_{t}^{0}(S_nu_0)\|_{B^s_{2,r}}
+\|\mathbf{S}_{t}^{0}(S_nu_0)-\mathbf{S}_{t}^{0}(u_0)\|_{B_{2,r}^{s}}\\
&\leq C\|(\mathrm{Id}-S_n)u_0\|_{B_{2,r}^{s}}+C2^{n}\alpha^{\fr12}.
\end{align*}
Choosing large $n$ and small $\alpha$ enables us to complete the proof of Theorem \ref{th1}.

\section{Proof of Theorem \ref{th2}}

\quad In this section, we will show that the zero-filter limit for the Camassa-Holm equation does not converges uniformly with respect to the initial data in $B^s_{2,r}(\R)$. First, we introduce a proposition which is useful in the proof of this theorem.

\begin{proposition}\label{pro 3.1}
Let Let $s>\frac32, 1\leq r<\infty$ and $\alpha\in[0,1)$. Assume that $\|u_0\|_{B^s_{2,r}}\approx 1$. Let $\mathbf{S}_{t}^{\mathbf{\alpha}}(u_0)\in\mathcal{C}([0, T];  B_{2, r}^s)$ and $\mathbf{S}_{t}^{0}(u_0)\in\mathcal{C}([0, T];  B_{2, r}^s)$ be the smooth solutions of \eqref{alpha-c} and \eqref{b} with the initial data $u_0$ respectively. Then we have for $t\in[0,T]$
\bbal
\f\|\mathbf{S}^\alpha_{t}(u_0)-u_0-t\mathbf{E}_0(\alpha)\g\|_{B_{2,r}^{s}}\leq Ct^{2}\mathbf{F}_0,
\end{align*}
where
\bbal
&\mathbf{E}_0(\alpha,u_0):=-3u_0\pa_xu_0-\alpha^2\partial^3_x \left(1-\alpha^2 \partial_x^2\right)^{-1}u_0^2-\frac{\alpha^2}{2}\partial_x \left(1-\alpha^2 \partial_x^2\right)^{-1} (\partial_xu_0)^2,\\
&\mathbf{F}_0(\alpha,u_0):=\f(\alpha\|u_0\|_{B_{2,r}^{s+1}}+1\g)\f(\alpha\|u_0\|_{B_{2,r}^{s+1}}+\|u_0\|_{B_{2,r}^{s-1}}\|u_0\|_{B_{2,r}^{s+1}}\g)+ \alpha\f(\|u_0\|_{B_{2,r}^{s+1}}+\|u_0\|_{B_{2,r}^{s-1}}\|u_0\|_{B_{2,r}^{s+2}}\g).
\end{align*}
\end{proposition}
\begin{proof} For simplicity, we denote $u(t)=\mathbf{S}^\alpha_t(u_0)$.
For $t\in[0,T]$, using Lemma \ref{le-pro}, Lemma \ref{lemm3} and the fundamental theorem of calculus in the time variable, we obtain
\bal\label{u1}
\|u(t)-u_0\|_{B^s_{2,r}}
&\leq \int^t_0\|\pa_\tau u\|_{B^s_{2,r}} \dd\tau
\nonumber\\&\leq \int^t_0\f(3\|u \pa_xu\|_{B^s_{2,r}}+\alpha^2\f\|\partial^3_x \left(1-\alpha^2 \partial_x^2\right)^{-1}u^2\g\|_{B^s_{2,r}}\g)\dd\tau\nonumber\\&\quad+\int^t_0\frac{\alpha^2}{2}\f\|\partial_x \left(1-\alpha^2 \partial_x^2\right)^{-1} (\partial_xu)^2\g\|_{B^s_{2,r}} \dd\tau
\nonumber\\&\lesssim t\f(\f\|u \pa_xu\g\|_{L_t^\infty B^s_{2,r}}+\alpha\f\|(\partial_xu)^2\g\|_{L_t^\infty B^s_{2,r}}\g)\nonumber\\
&\lesssim t\f(\|u\|_{L_t^\infty B_{2,r}^{s-1}}\|u\|_{L_t^\infty B_{2,r}^{s+1}}+\alpha\f\|\partial_xu\g\|_{L_t^\infty B_{2,r}^{s-1}}\f\|\partial_xu\g\|_{L_t^\infty B^s_{2,r}}\g)
\nonumber\\&\lesssim t\f(\|u_0\|_{B_{2,r}^{s-1}}\|u_0\|_{B_{2,r}^{s+1}}+\alpha\|u_0\|_{B_{2,r}^{s+1}}\g).
\end{align}
Following the same procedure of estimates as above, we have
\bal\label{u2}
\|u(t)-u_0\|_{B_{2,r}^{s-1}}
&\leq \int^t_0\|\pa_\tau u\|_{B_{2,r}^{s-1}} \dd\tau\nonumber\\
&\leq \int^t_0\f(3\|u \pa_xu\|_{B_{2,r}^{s-1}}+\alpha^2\f\|\partial^3_x \left(1-\alpha^2 \partial_x^2\right)^{-1}u^2\g\|_{B_{2,r}^{s-1}}\g)\dd\tau\nonumber\\
&\quad+\int^t_0\frac{\alpha^2}{2}\f\|\partial_x \left(1-\alpha^2 \partial_x^2\right)^{-1} (\partial_xu)^2\g\|_{B_{2,r}^{s-1}}  \dd\tau\nonumber\\
&\lesssim t\f(\f\|u \pa_xu\g\|_{L_t^\infty B_{2,r}^{s-1}}+\alpha\f\|(\partial_xu)^2\g\|_{L_t^\infty B_{2,r}^{s-1}}\g)\nonumber\\
&\lesssim t\f(\|u\|_{L_t^\infty B_{2,r}^{s-1}}\|u\|_{L_t^\infty B^s_{2,r}}+\alpha\|u\|^2_{L_t^\infty B^s_{2,r}}\g)
\nonumber\\&\lesssim t\f(\|u_0\|_{B_{2,r}^{s-1}}+\alpha\g),
\end{align}
and
\bal\label{u3}
\|u(t)-u_0\|_{B_{2,r}^{s+1}}
&\leq \int^t_0\|\pa_\tau u\|_{B_{2,r}^{s+1}} \dd\tau
\nonumber\\&\leq \int^t_0\f(3\|u \pa_xu\|_{B_{2,r}^{s+1}}+\alpha^2\f\|\partial^3_x \left(1-\alpha^2 \partial_x^2\right)^{-1}u^2\g\|_{B_{2,r}^{s+1}}\g)\dd\tau\nonumber\\
&\quad+\int^t_0\frac{\alpha^2}{2}\f\|\partial_x \left(1-\alpha^2 \partial_x^2\right)^{-1} (\partial_xu)^2\g\|_{B_{2,r}^{s+1}}  \dd\tau\nonumber\\
&\lesssim t\f(\f\|u \pa_xu\g\|_{L_t^\infty B_{2,r}^{s+1}}+\f\|(\partial_xu)^2\g\|_{L_t^\infty B_{2,r}^{s}}\g)\nonumber\\
&\lesssim t\f(\|u\|_{L_t^\infty B_{2,r}^{s-1}}\|u\|_{L_t^\infty B_{2,r}^{s+2}}+\|u\|_{L_t^\infty B^s_{2,r}}\|u\|_{L_t^\infty B_{2,r}^{s+1}}\g)
\nonumber\\&\lesssim t\f(\|u_0\|_{B_{2,r}^{s-1}}\|u_0\|_{B_{2,r}^{s+2}}+\|u_0\|_{B_{2,r}^{s+1}}\g).
\end{align}
For $t\in[0,T]$, using Lemma \ref{le-pro}, Lemma \ref{lemm3} and the fundamental theorem of calculus in the time variable again, we obtain
\bbal
&\|u(t)-u_0-t\mathbf{E}_0(\alpha,u_0)\|_{B^s_{2,r}}
\leq \int^t_0\|\pa_\tau u-\mathbf{E}_0(\alpha,u_0)\|_{B^s_{2,r}} \dd\tau
\\&\leq \int^t_0\f(3\|u\pa_xu-u_0\pa_xu_0\|_{B^s_{2,r}} +\alpha^2\f\|\partial^3_x \left(1-\alpha^2 \partial_x^2\right)^{-1}\f(u^2-u_0^2\g)\g\|_{B_{2,r}^{s}}\g)\dd\tau\\
&\quad+ \int^t_0\frac{\alpha^2}{2}\f\|\partial_x \left(1-\alpha^2 \partial_x^2\right)^{-1} \f((\partial_xu)^2-(\partial_xu_0)^2\g)\g\|_{B^s_{2,r}} \dd\tau\\
&\les \int^t_0\f(\|u^2-u_0^2\|_{B_{2,r}^{s+1}} +\alpha\f\|(\partial_xu)^2-(\partial_xu_0)^2\g\|_{B^s_{2,r}}\g) \dd\tau
\\&\les \int^t_0\f(\f(\alpha\|u_0\|_{B_{2,r}^{s+1}}+1\g)\|u(\tau)-u_0\|_{B^s_{2,r}} +\|u(\tau)-u_0\|_{B_{2,r}^{s-1}} \|u_0\|_{B_{2,r}^{s+1}}+\alpha\|u(\tau)-u_0\|_{B_{2,r}^{s+1}}\g)\dd \tau
\\&\les t^2\mathbf{F}_0(\alpha,u_0),
\end{align*}
where we have used \eqref{u1}-\eqref{u3} in the last step.
Thus, we complete the proof of Proposition \ref{pro 3.1}.
\end{proof}

\begin{lemma}\label{le3} Let $s>\frac32, 1\leq r<\infty$ and $f_n,g_n$ be defined as Lemma \ref{le2}. Assume that $s\in\R$ and $\alpha_n=2^{-n}$. Then there exists two positive constants $c=c(\phi)$ and  $C=C(\phi)$ such that
\bbal
&\liminf_{n\rightarrow \infty}\f\|\alpha^2_n\partial^2_x \left(1-\alpha^2_n \partial_x^2\right)^{-1}(g_n\pa_xf_n)\g\|_{B^s_{2,r}}\geq c.
\end{align*}
\end{lemma}
\begin{proof} By the construction of $f_n$ and $g_n$, one has (for more details see \cite{lyz})
\bbal
\mathrm{supp}\ \widehat{g_n\pa_xf_n}\subset \f\{\xi\in\R: \ \frac{17}{12}2^n-1\leq |\xi|\leq \frac{17}{12}2^n+1\g\}.
\end{align*}
According to Plancherel's identity, we obtain
\bbal
\f\|\alpha^2_n\partial^2_x \left(1-\alpha^2_n \partial_x^2\right)^{-1}(g_n\pa_xf_n)\g\|_{B^s_{2,r}}\approx\|g_n\pa_xf_n\|_{B^s_{2,r}}
\end{align*}
Using Lemma \ref{le2}, we complete the proof of this lemma.
\end{proof}
{\bf Proof of Theorem \ref{th2}.} Letting $\alpha_n=2^{-n}$ and setting $u^n_0=f_n+g_n$, it is easy to show that
\bbal
\|f_n\|_{B_{2,r}^{s+kn}}\lesssim 2^{kn}\quad\text{and}\quad \|g_n\|_{B_{2,r}^{s+kn}}\lesssim 2^{-n}\quad\text{for}\quad k\in\{-1,0,1,2\},
\end{align*}
which implies
\bal\label{chen-l-l}
\|u^n_0\|_{B_{2,r}^{s+kn}}\lesssim 2^{kn}.
\end{align}
Using \eqref{chen-l-l}, we can show that
\bbal
\mathbf{F}_0(\alpha_n,u^n_0)\lesssim  1.
\end{align*}
Notice that
\bbal
&\mathbf{S}^{\alpha_n}_{t}(u^n_0)=u^n_0+\underbrace{\mathbf{S}^{\alpha_n}_{t}(u^n_0)-u_0-t\mathbf{E}_0(\alpha_n,u^n_0)}_{=:\,\mathbf{I}_1}+t\mathbf{E}_0(\alpha_n,u^n_0),\\
&\mathbf{S}^{0}_{t}(u^n_0)=u^n_0+\underbrace{\mathbf{S}^{0}_{t}(u^n_0)-u_0-t\mathbf{E}_0(0,u^n_0)}_{=:\,\mathbf{I}_2}+t\mathbf{E}_0(0,u^n_0),\\
&\mathbf{E}_0(\alpha_n,u^n_0)-\mathbf{E}_0(0,u^n_0)=\alpha^2_n\partial^2_x \left(1-\alpha^2_n \partial_x^2\right)^{-1}(u^n_0\pa_xu^n_0)-\frac{\alpha^2_n}{2}\partial_x \left(1-\alpha^2_n \partial_x^2\right)^{-1} (\partial_xu^n_0)^2,
\end{align*}
and
$
u^n_{0}\pa_xu^n_{0}=g_n\pa_xf_n+f_n\pa_xf_n+u^n_{0}\pa_xg_n
$
, we deduce from Proposition \ref{pro 3.1}  that
\bal\label{yyh}
&\f\|\mathbf{S}^{\alpha_n}_{t}(u^n_0)-\mathbf{S}^0_{t}(u^n_0)\g\|_{B^s_{2,r}}
\geq~t\f\|\mathbf{E}_0(\alpha_n,u^n_0)-\mathbf{E}_0(0,u^n_0)\g\|_{B^s_{2,r}}-\f\|\mathbf{I}_1\g\|_{B^s_{2,r}}-\f\|\mathbf{I}_2\g\|_{B^s_{2,r}}\nonumber\\
\geq&~ t\f(\f\|\alpha^2_n\partial^2_x \left(1-\alpha^2_n \partial_x^2\right)^{-1}(u^n_0\pa_xu^n_0)\g\|_{B^s_{2,r}}
-\frac{\alpha^2_n}{2}\f\|\partial_x \left(1-\alpha^2_n \partial_x^2\right)^{-1} (\partial_xu_0)^2\g\|_{B^s_{2,r}}\g)-Ct^{2}\nonumber\\
\gtrsim&~ t\f(\f\|\alpha^2_n\partial^2_x \left(1-\alpha^2_n \partial_x^2\right)^{-1}(g_n\pa_xf_n)\g\|_{B^s_{2,r}}
-\|f_n\pa_xf_n\|_{B^s_{2,r}}-\|u^n_{0}\pa_xg_n\|_{B^s_{2,r}}-2^{-n}\f\|(\partial_xu^n_0)^2\g\|_{B^s_{2,r}}\g)-t^{2}.
\end{align}
Using Lemma \ref{le-pro} and Lemma \ref{le2}, we get
\bbal
&\|f_n\pa_xf_n\|_{B^s_{2,r}}\leq C\|f_n\|_{L^\infty}\|f_n\|_{B_{2,r}^{s+1}}+C\|\pa_xf_n\|_{L^\infty}\|f_n\|_{B_{2,r}^{s}}\leq C2^{-n(s-1)},\\
&\f\|u^n_{0}\pa_xg_n\g\|_{B^s_{2,r}}\leq C\|u^n_0\|_{B^s_{2,r}}\|g_n\|_{B_{2,r}^{s+1}}\leq C2^{-n},\\
&\|(\partial_xu^n_0)^2\|_{B^s_{2,r}}\leq C\|\partial_xu^n_0\|_{L^\infty}\|u^n_0\|_{B_{2,r}^{s+1}}\leq 2^n(2^{-n}+2^{-n(s-1)})\leq C(1+2^{-n(s-2)}).
\end{align*}
Plugging the above estimates into \eqref{yyh} and combining Lemma \ref{le3} yields that
\bbal
\liminf_{n\rightarrow \infty}\|\mathbf{S}^{\alpha_n}_{t}(u^n_0)-\mathbf{S}^0_{t}(u^n_0)\|_{B^s_{2,r}}\gtrsim t\quad\text{for} \ t \ \text{small enough}.
\end{align*}
This completes the proof of Theorem \ref{th2}.

\section*{Acknowledgments}
M. Li is supported by the Jiangxi Provincial Natural Science Foundation (20232BAB201013). J. Li is supported by the National Natural Science Foundation of China (12161004), Training Program for Academic and Technical Leaders of Major Disciplines in Ganpo Juncai Support Program
(20232BCJ23009) and Jiangxi Provincial Natural Science Foundation (20224BAB201008).

\section*{Conflict of interest}
The authors declare that they have no conflict of interest.

\section*{Data Availability} Data sharing is not applicable to this article as no new data were created or analyzed in this study.

\end{document}